\documentclass[12pt,twoside]{article}
\title{Rings over which all modules are  strongly Gorenstein
projective}
\date{}
\usepackage{amsfonts}
\usepackage{amsmath}
\usepackage{amssymb}
\usepackage{latexsym}
\usepackage{graphicx}
\newtheorem{thm}{\bf Theorem}[section]
\newtheorem{cor}[thm]{\bf Corollary}
\newtheorem{lem}[thm]{\bf Lemma}
\newtheorem{prop}[thm]{\bf Proposition}

\newtheorem{defns}[thm]{\bf Definitions}

\catcode`\ç=13
\defç{\c{c}}
\catcode`\é=13
\defé{\'e}
\catcode`\à=13
\defà{\`a}
\catcode`\è=13
\defè{\`e}
\catcode`\â=13
\defâ{\^a}
\catcode`\ù=13
\defù{\`u}
\catcode`\ê=13
\defê{\^e}
\catcode`\î=13
\defî{\^\i}
\catcode`\ô=13
\defô{\^o}
\newcommand{\field}[1]{\mathbb{#1}}

\newcommand{\Z }{\field{Z}}

\def\proof{{\parindent0pt {\bf Proof.\ }}}

\def\Im{{\rm Im}}

\def\Ker{{\rm Ker}}

\def\Ext{{\rm Ext}}

\def\Hom{{\rm Hom}}

\newcommand{\cqfd}
{\hspace{1cm}
\rule{2mm}{2mm}%
\medbreak%
\par%
}
\def\1{{\noindent\rm (1)}}
\def\2{{\noindent\rm (2)}}
\def\3{{\noindent\rm (3)}}
\def\4{{\noindent\rm (4)}}
\def\5{{\noindent\rm (5)}}
\begin{document}
\thispagestyle{empty}

%%%%%%%%%%%%%%%%%%%%%%%%%%%%%%%%%%%%%%%%%%%%%%%%%%%%%%%%%
%%%%%%%%%%%%%%%%%%%%%%%%%%%%%%%%%%%%%%%%%%%%%%%%%%%%%%%%%
%%%%%%%%%%%%%%%%%%%%%%%%%%%%%%%%%%%%%%%%%%%%%%%%%%%%%%%%%
%%%TITLE%%%%%%%%%%%%%%%%%%%%%%%%%%%%%%%%%%%%%%%%%%%%%%%%%
\maketitle \vspace*{-1.5cm}
\begin{center}{\large\bf Driss Bennis, Najib Mahdou, and Khalid Ouarghi}
%%%%%%%%%%%%%%%%%%%%%%%%%%%%%%%%%%%%%%%%%%%%%%%%%%%%%%%%%
%%%%%%%%%%%%%%%%%%%%%%%%%%%%%%%%%%%%%%%%%%%%%%%%%%%%%%%%%
%%%%%%%%%%%%%%%%%%%%%%%%%%%%%%%%%%%%%%%%%%%%%%%%%%%%%%%%%
%%%NAMES%%%%%%%%%%%%%%%%%%%%%%%%%%%%%%%%%%%%%%%%%%%%%%%%%

\bigskip

%%%%%%%%%%%%%%%%%%%%%%%%%%%%%%%%%%%%%%%%%%%%%%%%%%%%%%%%%
%%%%%%%%%%%%%%%%%%%%%%%%%%%%%%%%%%%%%%%%%%%%%%%%%%%%%%%%%
%%%%%%%%%%%%%%%%%%%%%%%%%%%%%%%%%%%%%%%%%%%%%%%%%%%%%%%%%
%%%%%%%%%%%%ADDRESSES%%%%%%%%%%%%%%%%%%%%%%%%%%%%%%%%%%%%%%%%%%%%%
 \small{Department of Mathematics, Faculty of Science and Technology of Fez,\\ Box 2202, University S. M.
Ben Abdellah Fez, Morocco\\[0.12cm]
driss\_bennis@hotmail.com\\
mahdou@hotmail.com\\
ouarghi.khalid@hotmail.fr }
\end{center}

\bigskip\bigskip
%%%%%%%%%%%%%%%%%%%%%%%%%%%%%%%%%%%%%%%%%%%%%%%%%%%%%%%%%

%%%%%%%%%%%%%%%%%%%%%%%%%%%%%%%%%%%%%%%%%%%%%%%%%%%%%%%%
%%%%%%%%%%%%%%%%%%%%%%%%%%%%%%%%%%%%%%%%%%%%%%%%%%%%%%%%%
%%%ABSTRACT%%%%%%%%%%%%%%%%%%%%%%%%%%%%%%%%%%%%%%%%%%%%%%
\noindent{\large\bf Abstract.}  One of the main results of this
paper is the characterization of  the rings over which all modules
are strongly Gorenstein projective. We show that these kinds of
rings are very particular cases of the well-known quasi-Frobenius
rings. We give examples of  rings over which all modules are
Gorenstein projective but not necessarily strongly Gorenstein
projective.\bigskip

%%%%%%%%%%%%%%%%%%%%%%%%%%%%%%%%%%%%%%%%%%%%%%%%%%%%%%%%%
\small{\noindent{\bf Key Words.} (Strongly) Gorenstein projective,
injective, and flat modules; quasi-Frobenius rings.}
\bigskip\bigskip
 %%%%%%%%%%%%%%%%%%%%%%%%%%%%%%%%%%%%%%

%%%%%%%%%%%%%%%%%%%%%%%%%%%%%%%%%%%%%%%%%%%%%%%%%%%%%%%%
%%%%%%%%%%%%%%%%%%%%%%%%%%%%%%%%%%%%%%%%%%%%%%%%%%%%%%%%
%%%%%%%%%%%%%%%%%%%%%%%%%%%%%%%%%%%%%%%%%%%%%%%%%%%%%%%%

%%%%%%%%%%%%%%%%%%%%%%%%%%%%%%%%%%%%%%%%%%%%%%%%%%%%%%%%%
%%%INTRODUCTION%%%%%%%%%%%%%%%%%%%%%%%%%%%%%%%%%%%%%%%%%%
\begin{section}{Introduction} Throughout this paper all rings are
commutative with identity element and all modules are unital. It
is convenient to use ``$m$-local" to refer to (not necessarily
Noetherian) rings with a unique maximal ideal $m$.\bigskip

For background on the following definitions, we refer the reader
to \cite{BM, LW, Rel-hom, HH}.

\begin{defns}\label{Def-Goren-Mod}
  \textnormal{A module $M$ is said to be
\textit{Gorenstein projective}, if there exists an exact sequence
of  projective modules $$\mathbf{P}=\ \cdots\rightarrow
P_1\rightarrow P_0 \rightarrow P^0
         \rightarrow P^1 \rightarrow\cdots$$ such that  $M \cong \Im(P_0
\rightarrow \nolinebreak P^0)$ and such that $\Hom( -, Q) $ leaves
the sequence $\mathbf{P}$ exact whenever $Q$ is a projective
module.\\
\indent The exact sequence $\mathbf{P}$ is called a \textit{complete
projective} resolution.\\
\indent The \textit{Gorenstein injective}  modules are defined
dually.}
\end{defns}

Recently in \cite{BM}, the authors studied a simple particular
case of Gorenstein projective and injective modules, which are
defined, respectively, as follows:

\begin{defns}[\cite{BM}]\label{defSG}
 \textnormal{A module $M$ is said to be \textit{strongly
Gorenstein projective}, if there exists a complete projective
resolution of the form $$ \mathbf{P}=\
\cdots\stackrel{f}{\longrightarrow}P\stackrel{f}{\longrightarrow}P\stackrel{f}{\longrightarrow}P
\stackrel{f}{\longrightarrow}\cdots $$ such that  $M \cong
\Im(f)$.\\
\indent  The exact sequence $\mathbf{P}$ is called a \textit{strongly
complete
projective} resolution.\\
\indent  The \textit{strongly Gorenstein injective} modules are
defined dually.}
\end{defns}

The principal role of the strongly Gorenstein projective and
injective modules is to give a simple characterization of
Gorenstein projective  and injective modules, respectively, as
follows:

\begin{thm}[\cite{BM}, Theorem 2.7] \label{thm-car-G-SG}
A module is Gorenstein projective  (resp., injective) if and only
if it is a direct summand of a strongly Gorenstein projective
(resp., injective) module.
\end{thm}

The important of this last result manifests in showing that  the
strongly Gorenstein projective and  injective modules have simpler
characterizations than their Gorenstein correspondent modules. For
instance:

\begin{prop}[\cite{BM}, Proposition  2.9] \label{pro-cara-SG-pro}
A module $M$ is strongly Gorenstein projective if and only if
there exists a short exact sequence of modules $$0\rightarrow
M\rightarrow P\rightarrow M\rightarrow 0,$$ where $P$ is
projective, and $\Ext(M,Q)=0$ for any   projective module $Q$.
\end{prop}

The aim of this paper is to investigate the two following classes
of rings:
\begin{enumerate}
    \item The rings over which all modules are Gorenstein
    projective (resp., injective), which are called \textit{G-semisimple} rings (please see Proposition \ref{pro-equi-1}).
    \item The rings over which all modules are strongly Gorenstein
    projective (resp., injective), which are called \textit{SG-semisimple} rings (please see Proposition \ref{pro-equi-2}).
\end{enumerate}

In Section 2, we show that the  G-semisimple rings are just the
well-known quasi-Frobenius rings; i.e.,  Noetherian and
self-injective rings. The SG-semisimple rings are then particular
cases of the quasi-Frobenius rings. In Section 3, we characterize
the SG-semisimple rings. Namely, we show that an $m$-local ring is
SG-semisimple if and only if it has at most one proper nonzero
ideal; in general, a ring is SG-semisimple if and only if it is a
finite direct product of local SG-semisimple rings.\bigskip

Before starting, we need to recall some useful results about
quasi-Frobenius rings (for more details about this kind of rings
see for example \cite{AF} and \cite{NY}). The quotient ring
$R/I$, where $R$ is a principal ideal domain and $I$ is any
nonzero ideal of $R$, is a classical example of a quasi-Frobenius
ring \cite[Exercise 9.24]{Rot}. The quasi-Frobenius rings have
several characterizations. Here, we need the following:

\begin{thm}[\cite{NY}, Theorems 1.50,
7.55, and 7.56] \label{thm-car-QF} For a ring $R$, the following
are equivalent:
\begin{enumerate}
    \item $R$ is quasi-Frobenius;
    \item $R$ is Artinian and self-injective;
    \item Every projective $R$-module is injective;
    \item Every injective $R$-module is projective;
    \item $R$ is Noetherian and, for every ideal $I$, $Ann (Ann(I))=I$,
     where $Ann(I)$ denotes the annihilator of $I$.
\end{enumerate}
\end{thm}

The quasi-Frobenius rings are particular cases of the perfect
rings; i.e., the rings over which all flat modules are projective.
Namely, a ring is quasi-Frobenius if and only if it is perfect and
self-injective \cite[Theorem 6.39]{NY}. The perfect rings are
introduced by Bass in \cite{bass}. They have the following
characterizations (needed later):

\begin{thm}[\cite{bass}, Theorem P and Example 6, p. 476] \label{thm-car-perfect}
For a ring $R$, the following are equivalent:
\begin{enumerate}
    \item $R$ is perfect;
    \item Every direct limit (with directed index set) of projective $R$-modules is projective;
    \item $R$ is a finite direct product of local rings, each with
$T$-nilpotent maximal ideal (i.e., if we pick a sequence $a_1,
a_2,...$ of elements in the maximal ideal, then for some index
$j$, $a_1a_2...a_j=0$).
\end{enumerate}
\end{thm}

From Theorems \ref{thm-car-QF} and \ref{thm-car-perfect} above and
\cite[Lemma 5.64]{NY}, we may give the following structural
characterization of quasi-Frobenius rings, which will be used
later:

\begin{prop}\label{thm-car-QF-product}
A ring $R$ is quasi-Frobenius if and only if
$R=R_1\times\cdots\times R_n$, where each $R_i$ is a local
quasi-Frobenius ring.
\end{prop}

\end{section}
%%%%%%%%%%%%%%%%%%%%%%%%%%%%%%%%%%%%%%%%%%%%%%%%%%%%%%%%%%%%%%%%%%%%%%%%%%%%
%%%%%%%%%%%%%%%%%%%%%%%%%%%%%%%%%%%%%%%%
%%%%%%%%%%%%%%%%%%%%%%%%%%%%%%%%
%%%%%%%%%%%%%%%%%%%%%%%%                  $$$$                      $$$$
%%%%%%%%%%%%%%%%%%%                             S 2  G--SYZYGY THEOREM
%%%%%%%%%%%%%%%%%%%%%%%%                  $$$$                      $$$$
%%%%%%%%%%%%%%%%%%%%%%%%%%%%%%%%
%%%%%%%%%%%%%%%%%%%%%%%%%%%%%%%%%%%%%%%%
%%%%%%%%%%%%%%%%%%%%%%%%%%%%%%%%%%%%%%%%%%%%%%%%%%%%%%%%%%%%%%%%%%%%%%%%%%%%
%%%%%%%%%%%%%%%%%%%%%%%%%%%%%%%%%%%%%%%%%%%%%%%%%%%%%%%%%%%%
\begin{section}{G-semisimple rings}
In this section we investigate the G-semisimple rings; i.e., the
rings that satisfy each of the  following equivalent conditions:

\begin{prop}\label{pro-equi-1}
Let  $R $ be a ring. The following are equivalent:
\begin{enumerate}
    \item Every $R$-module is Gorenstein projective;
    \item Every $R$-module is Gorenstein injective.
\end{enumerate}
\end{prop}
\proof We prove the implication $(1) \Rightarrow (2)$, and
the proof of the converse implication is analogous.\\
Assume that every module is Gorenstein projective. Then, any
injective module is projective (since, as a Gorenstein projective
module, it embeds in a projective module). This is equivalent, by
Theorem \ref{thm-car-QF}, to say that every projective module is
injective. Then, every complete projective resolution is also a
complete injective resolution, and therefore, every $R$-module is
Gorenstein injective.\cqfd\bigskip

Note that the equivalence in Proposition \ref{pro-equi-1} is
already known when $R$ is  Noetherian, and that each of the
conditions (1) and (2) is equivalent to  the ring being
quasi-Frobenius (see for example \cite[Theorem 12.3.1]{Rel-hom}).
Next, we show how  Proposition \ref{pro-equi-1} and its proof show
that a G-semisimple ring is the same as a quasi-Frobenius ring.

\begin{thm}\label{thm-car-G-semi}
For any ring $R$, the following are equivalent:
\begin{enumerate}
  \item  $R$ is G-semisimple;
  \item  Every Gorenstein injective $R$-module is Gorenstein projective;
  \item  Every strongly Gorenstein injective $R$-module is strongly Gorenstein projective;
  \item  Every Gorenstein projective $R$-module is Gorenstein injective;
  \item  Every strongly Gorenstein injective $R$-module is strongly Gorenstein projective;
  \item  $R$ is quasi-Frobenius.
\end{enumerate}
\end{thm}
\proof  First Note that a G-semisimple ring is Noetherian. Indeed,
from the proof of  Proposition \ref{pro-equi-1}, we have that if
$R$ is a G-semisimple ring, then every projective $R$-module is
injective. This means from Theorem \ref{thm-car-QF}, that $R$ is
quasi-Frobenius and so is Noetherian. This gives a proof of the
implication $(1)\Rightarrow (6)$. For the proof of the remains
implications use also Proposition \ref{pro-equi-1} and its
proof.\cqfd\bigskip

We have the following relationship between  semisimple rings and
the G-semisimple rings; compare to \cite[Exercise 9.2]{Rot}.

\begin{prop}
A G-semisimple ring is semisimple if and only if it has finite
global dimension.
\end{prop}
\proof Follows from  the  fact that a Gorenstein projective
 module is projective  if and
only if it has finite projective  dimension \cite[Proposition
2.27]{HH}.\cqfd\bigskip

 Finally, it is important to say that there exist numerous  examples of
G-semisimple rings which are not semisimple, for instance
$\Z/4\Z$.

\end{section}
%%%%%%%%%%%%%%%%%%%%%%%%%%%%%%%%%%%%%%%%%%%%%%%%%%%%%%%%%%%%
%%%%%%%%%%%%%%%%%%%%%%%%%%%%%%%%%%%%%%%%%%%%%%%%%%%%%%%
%%%%%%%%%%%%%%%%%%
%%%%%%%%%%%%%%%%%%%%%%   3.      SG-semisimple rings
%%%%%%%%%%%%%%%
%%%%%%%%%%%%%%%%%%%%%%%%%%%%%%%%%%%%%%%%%%%%%%%%%%%%%%%%%
%%%%%%%%%%%%%%%%%%%%%%%%%%%%%%%%%%%%%%%%%%%%%%%%%%%%%%%%%
\begin{section}{SG-semisimple rings}   We investigate, in this
section, the SG-semisimple rings; i.e., rings that satisfy each of
the following equivalent conditions.

\begin{prop}\label{pro-equi-2}
Let  $R $ be a ring. The following are equivalent:
\begin{enumerate}
    \item Every $R$-module is strongly Gorenstein projective;
    \item Every $R$-module is strongly Gorenstein injective.
\end{enumerate}
\end{prop}
\proof It suffices to prove the implication $(1) \Rightarrow
(2)$, and the proof of the converse implication is analogous.\\
Assume that every module is strongly Gorenstein projective. Then,
by Theorem  \ref{thm-car-G-semi}, $R$ is G-semisimple (i.e.,
quasi-Frobenius). Thus,  we can show that a strongly complete
projective resolution is also a strongly complete injective
resolution.\cqfd\bigskip

Naturally, an SG-semisimple ring is G-semisimple (i.e.,
quasi-Frobenius). Later, we give examples of SG-semisimple rings
and other examples of G-semisimple rings which are not
SG-semisimple (see Corollaries \ref{exm-1} and \ref{exm-2}).
Before that, we give a characterization of SG-semisimple rings. We
begin by a structure theorem. For that, we need the following
lemma.

\begin{lem}\label{lem-cara-sG-proj-dir-prod}
Let $R=R_1\times\cdots\times R_n$ be a finite direct product of
rings $R_i$. An $R$-module $M$ is (strongly) Gorenstein projective
if and only if $M=M_1\oplus\cdots\oplus M_n$, where each $M_i$ is
a (strongly) Gorenstein projective $R_i$-module.
\end{lem}
\proof This follows from the structure of (projective) modules and
homomorphisms over a finite direct product of rings (see for
example \cite[Subsection 2.6]{BK}).\cqfd

\begin{thm}\label{thm-SG-semi-1}
A ring $R$ is SG-semisimple if and only if
$R=R_1\times\cdots\times R_n$, where each $R_i$ is a local
SG-semisimple ring.
\end{thm}
\proof The result is a  consequence of Proposition
\ref{thm-car-QF-product} and Lemma \ref{lem-cara-sG-proj-dir-prod}
above.\cqfd\bigskip

Theorem \ref{thm-SG-semi-1} leads us to restrict the study of the
SG-semisimple rings to the local SG-semisimple rings.

\begin{lem}\label{lem-1}
Let $R$  be an $m$-local ring and let $x\neq 0$ be a zero-divisor
element of $R$. If the ideal $xR$ is strongly Gorenstein
projective, then $Ann(xR)\cong xR$ and therefore
$Ann(Ann(xR))=Ann(xR)$.\\
Particularly, if $xR=m$, we get $Ann(m)=m$.
\end{lem}
\proof Since $xR$ is strongly Gorenstein projective, there exists,
by \cite[Proposition 2.9]{BM}, a short exact sequence of
$R$-modules $$(\star)\qquad\qquad 0\rightarrow xR \rightarrow P
\rightarrow xR \rightarrow  0,$$ where $P$ is projective, then
free (since $R$ is $m$-local). In the sequence $(\star)$ $xR$ is
finitely generated, then so is the free $R$-module $P$. Thus,
there exists a nonzero positive integer $n$ such that $P\cong
R^n$. Hence, we get the following exact sequence:
$$(\star\star)\qquad\qquad 0\rightarrow xR \rightarrow R^n \rightarrow
xR \rightarrow  0.$$
 Consider also the following canonical short
exact sequence of $R$-modules:
$$ 0\rightarrow Ann(xR) \rightarrow R \rightarrow
xR \rightarrow  0.$$ From Schanuel's lemma \cite[Theorem
3.62]{Rot}, we have: $$Ann(xR) \oplus R^{n}\cong R\oplus( xR).$$
Then, since $R$ is $m$-local, the  minimal generating sets of both
$Ann(xR) \oplus R^{n}$ and $ R\oplus( xR) $ have the same numbers
of elements which is necessarily $2$.  On the other hand, since
$x$ is a zero-divisor element of $R$, $Ann(xR)\neq 0$. Thus,
$Ann(xR)$ is generated by at least one element, and so $Ann(xR)
\oplus R^{n}$ is generated by at least $n+1$ elements. Then, by
the reason above, $n$ must equal $1$. So the sequence
$(\star\star)$ becomes:
$$0\rightarrow xR \rightarrow R \stackrel{f}\rightarrow
xR \rightarrow  0.$$ Now, let $\alpha \in R$ with $f(1)=\alpha x$.
Since $f$ is surjective, there exists $\beta\in R$ such that
$f(\beta)=x$. So, $x=\beta\alpha x$, and then $(1-
\beta\alpha)x=0$, which means that $(1- \beta\alpha)\in
Ann(xR)\subseteq m$. Then, $\beta\alpha$ is invertible and so is
$\alpha$. This implies that:
$$\Ker\,f=\{y\in R\,|\,0=yf(1)=y \alpha x\}= Ann(xR).$$
 Consequently, $ xR \cong \Ker\,f= Ann(xR)$. Therefore,
$Ann(Ann(xR))=Ann(xR)$, as desired.\\
\indent Now, if $m=xR$, then $m=xR\subseteq
Ann(Ann(xR))=Ann(xR)\subseteq m$.\cqfd

\begin{lem}\label{lem-2}
Let $R$  be an $m$-local ring and let $I$ be a nonzero proper
ideal of $R$. If $R/I$ is a strongly Gorenstein projective
$R$-module, then $I$ is a cyclic strongly Gorenstein projective
ideal generated by a zero-divisor element of $R$.
\end{lem}
\proof Since $R$ is a $m$-local ring and similarly to the first
part of the proof of Lemma \ref{lem-1} above, we get a  short
exact sequence of $R$-modules:
$$(\ast)\qquad\qquad 0\rightarrow R/I \rightarrow R^n \rightarrow
R/I \rightarrow  0,$$ where $n$ is a nonzero positive integer. \\
And also, the same argument as in the proof of Lemma \ref{lem-1}
above, and using the short exact sequence of $R$-modules:
$$ (\ast\ast)\qquad\qquad 0\rightarrow I \rightarrow R \rightarrow
R/I \rightarrow  0,$$ we get $n=1$ and $I=xR$ for some
zero-divisor element $x$ of $R$.\\
Now, to show that $I$ is strongly Gorenstein projective, note at
first that it is Gorenstein projective (by the sequence $
(\ast\ast)$ and from \cite[Theorem 2.5]{HH}). Then, $ \Ext
(I,P)=0$ for every projective $R$-module $P$ (by \cite[Proposition
2.3]{HH}). On the other hand,   the two sequences $(\ast)$ and $
(\ast\ast)$ with the Horseshoe Lemma  \cite[Lemma 6.20]{Rot} give
the following commutative diagram with exact columns and rows:
$$\begin{array}{ccccccccc}
    &  & 0 &  & 0&  &0 &  &  \\
    &  & \downarrow &  &\downarrow&  &\downarrow &  &  \\
   0 & \rightarrow & I&\rightarrow  &Q& \rightarrow &I&\rightarrow  & 0 \\
    &  & \downarrow &  &\downarrow&  &\downarrow &  &  \\
   0 & \rightarrow & R&\rightarrow  &R\oplus R& \rightarrow &R &  \rightarrow&  0\\
    &  & \downarrow &  &\downarrow&  &\downarrow &  &  \\
   0 & \rightarrow & R/I &\rightarrow  &R& \rightarrow &R/I & \rightarrow & 0  \\
    &  & \downarrow &  &\downarrow&  &\downarrow &  &  \\
    &  & 0 &  & 0&  &0 &  &
\end{array}$$
Note that $Q$ is a projective $R$-module. Therefore, by the top
horizontal sequence and Proposition \ref{pro-cara-SG-pro}, $I$ is
a strongly Gorenstein projective ideal.\cqfd

\begin{lem}\label{lem-3}
If $R$  is a  local G-semisimple ring, then every $R$-module $M$
is of the form: $M=R^{(I)}\oplus N$, where $I$ is an index set and
$N$ is an $R$-module with $Ann(x)\not = 0$ for every element $x$
of $ N$.
\end{lem}
\proof We may assume that $M$ admits an element $x$ such that
$rx\not =0$ for all $0\neq r\in R$. Consider the set $E$ of  all
free submodules of $M$. The set $E $ is not empty, since $xR$ is a
free submodule of $M$. On the other hand, since $R$ is a  local
G-semisimple ring and from Theorem \ref{thm-car-perfect}, a direct
limit of free $R$-modules is a free $R$-module. Then, for every
subchain $E_i$ of $E$, $\cup E_i$ is a free submodule of $M$.
Then, by Zorn's lemma, $E$ admits a maximal element $F$. We may
set $F=R^{(I)}$ which is injective (since $R$ is G-semisimple).
Then, $F$ is a direct summand of $M$ and so $M=F\oplus N$ for some
$R$-module $N$. If there exists $x\in N$ such that $rx\not =0$ for
all $r\in R$, then $xR\cong R$ is injective and then a direct
summand of $N$. Hence, there exists an $R$-module $N'$ such that
$N=xR\oplus N'$, and so $M=F\oplus N= F\oplus xR\oplus N'$. But,
the free submodule $F\oplus xR$ of $M$ contradicts the maximality
of $F$.\cqfd\bigskip

The main result, in this section, is the following
characterization of local SG-semisimple rings.

\begin{thm}\label{thm-cara-local-SG-semi}
Let $R$  be an $m$-local ring. The following are equivalent:
\begin{enumerate}
    \item $R$ is SG-semisimple;
    \item $R/m$ is a strongly Gorenstein projective $R$-module;
    \item $R$ has a most one nonzero proper ideal (which is necessarily $m$).
\end{enumerate}
\end{thm}
\proof
 $(1)\Rightarrow (2)$. By definition.\\
\indent $(2)\Rightarrow (3)$. From Lemma \ref{lem-2}, $m=xR$ is a
cyclic strongly Gorenstein projective ideal and $x$ is
zero-divisor. Then, by Lemma \ref{lem-1}, $m^2=0$. Therefore, a
standard argument shows that either $m=0$ or $m$ is the unique nonzero proper ideal of $R$.\\
\indent $(3)\Rightarrow (1)$. We may assume that $R$ is not a
field. Clearly $m=xR$ (for some $ 0 \neq x \in R$) and $m^2=0$.
Then, from Theorem \ref{thm-car-QF}, $R$ is G-semisimple (i.e.,
quasi-Frobenius), and so  $m$ is a Gorenstein projective  ideal of
$R$. Hence, by \cite[Proposition 2.3]{HH}, $\Ext(m,Q)=0$ for every
projective $R$-module $Q$.  Then, by the short exact sequence
$$0\rightarrow Ann(m)=m \rightarrow R \rightarrow m \rightarrow
0,$$ and from Proposition \ref{pro-cara-SG-pro}, $m$ is a
strongly Gorenstein projective $R$-module.\\
\indent Now, consider an arbitrary $R$-module $M$. By Lemma
\ref{lem-3}, there exists an index set $I$ such that $M\cong
R^{(I)}\oplus N$, where $N$ is an $R$-module  with $Ann(y)\not =
0$ for every nonzero element $y\in N$. Then, necessarily $xN=0$,
and so $N\cong (R/m)^{(J)}$ for some index set $J$. Since
$R/m\cong m$ is a strongly Gorenstein projective $R$-module  and,
by \cite[Proposition 2.2]{BM}, $N$ is a strongly Gorenstein
projective $R$-module. Therefore, $M$ is a strongly Gorenstein
projective $R$-module.\cqfd

\begin{cor}\label{cor-thm1-thm2}
A ring $R$ is SG-semisimple if and only if
$R=R_1\times\cdots\times R_n$, where each $R_i$ is a ring with at
most one nonzero proper ideal.
\end{cor}
\proof Combine Theorems \ref{thm-SG-semi-1} and
\ref{thm-cara-local-SG-semi}.\cqfd\bigskip

We end with some examples of G-semisimple and SG-semisimple rings.

\begin{cor}\label{exm-1}
For every principal ideal domain $R$ and every nonzero prime ideal
$p$ of $R$, the ring $R/p^2$ is  local SG-semisimple.
\end{cor}

The following result shows how to construct G-semisimple rings
which are not SG-semisimple.

\begin{cor}\label{exm-2}
For every principal ideal domain $R$ and every nonzero prime ideal
$p$ of $R$,   the ring $R/p^n$, where $n\geq 3$, is a local
G-semisimple, but it is not SG-semisimple.
\end{cor}

%%%%%%%%%%%%%%%%%%%%%%%%%%%%%%%%%%%%%%%%%%%%%%%%%%%%%%%%
\end{section}
%%%%%%%%%%%%%%%%%%%%%%%%%%%%%%%%%%%%%%%%%%%
%%%%%%%%%%%%%%%%%%%%%%%%%%%%%%%%%%%%%%%%%%%%%%%%%%%%%%%%
%%%%%%%%%%%%%%%%%%%%%%%%%%%%%%%%%%%%%%%%%%%%%%%%%%%%
%%%%%%%%%%%%%%%%%%%%%%%%%%%%%%%%%%%%%%%%%%%%%%%%%%%%
%%%%%%%%%%%%%%%%%%%%%%%%%%%%%%%%%%%%%%%%%%%%%%%%%%%%%%%%%%%%
%%%%%%%%%%%%%%%%%%%%%%%%%%%%%%%%%%%%%%%%%%%%%%%%%%%%%%%%%
%%%%%%%%%%%%%%%%%%%%%%%%%%%%%%%%%%%%%%%%%%%%%%%%%%%%%%%%%
%%%REFERENCES%%%%%%%%%%%%%%%%%%%%%%%%%%%%%%%%%%%%%%%%%%%%
%%%%%%%%%%%%%%%%%%%%%%%%%%%%%%%%%%%%%%%%%%%%%%%%%%%%%%%

%%%%%%%%%%%%%%%%%%%%%%%%%%%%%%%%%%%%%%%%%%%%%%%%%%%%%%%%
\end{document}